\documentclass[12pt,twoside]{article}

\usepackage{amssymb}
\usepackage{amsmath,amsbsy,amsfonts}
\usepackage{color}
\usepackage{graphicx}
\date{}
\def\nd{\noindent}
\def\thend{\rule{3mm}{3mm}}

\newtheorem{theorem}{Theorem}[section]
\newtheorem{prop}{Proposition}[section]
\newtheorem{lem}{Lemma}[section]
\newtheorem{rmk}{Remark}[section]

\newenvironment{dem}{\nd{\bf Proof. }}{\hskip.3cm\thend}

\title{Existence of ground state solutions to Dirac equations with vanishing potentials at infinity}

\author{Giovany M. Figueiredo\thanks{Partially supported by CNPq/PQ
301242/2011-9 and FAPESP 2015/12476-5}\\
\noindent Faculdade de Matem\'{atica}\\
\noindent Universidade Federal do Par\'a \\
\noindent 66075-110, Bel\'{e}m - Pa, Brazil.\\
\noindent e-mail: {\tt{giovany@ufpa.br}}\\
\mbox{}\\
and \\
Marcos T. O. Pimenta\thanks{Corresponding author } \thanks{Supported by FAPESP 2014/16136-1 and CNPq 442520/2014-0}\\
\noindent Departamento de Matem\'atica e Computa\c{c}\~ao\\
\noindent Faculdade de Ci\^encias e Tecnologia\\
\noindent Universidade Estadual Paulista - Unesp\\
\noindent 19060-900, Presidente Prudente - SP, Brazil.\\
\noindent e-mail: {\tt{pimenta@fct.unesp.br}}\\
}
\date{}

\begin{document}

\maketitle

\begin{abstract}
In this work we study the existence of ground-state solutions of Dirac equations with potentials which are allowed to vanish at infinity. The approach is based on minimization of the energy functional over a generalized Nehari set. Some conditions on the potentials are given in order to overcome the lack of compactness.

\end{abstract}



\section{Introduction}

\hspace{0.5cm} In 1928 Dirac proposed a model to the quantum mechanics which, in contrast to the Schr\"odinger theory, takes into account the Relativity Theory. More specifically, he proposed a model to describe the evolution of a free relativistic particle, given by
\begin{equation}
i\hbar\frac{\partial\psi}{\partial t} = D_c\psi,
\label{Diraclivre}
\end{equation}
where the operator $D_c$ is given by
$$D_c = -i c\hbar \boldsymbol{\alpha}\cdot\nabla + mc^2 \beta = -ic\hbar \sum_{k=1}^3 \alpha_k\partial_k + mc^2\beta$$
and $\boldsymbol{\alpha} = (\alpha_1,\alpha_2,\alpha_3)$, $\beta$  satisfy the anticomutation relations
\begin{equation}
\left\{
\begin{array}{rll}
\alpha_k\alpha_l + \alpha_l\alpha_k & = & 2\delta_{kl}I,\\
\alpha_k\beta + \beta\alpha_k & = & 0,\\
\beta^2 & = & I,
\end{array}\right.
\label{Pauli}
\end{equation}
where $I$ denotes the identity matrix. It can be proved that the least dimension where (\ref{Pauli}) can hold is $N = 4$, where $\alpha_i$ and $\beta$ are four-dimensional complex matrices given by
$$\beta = \left(
\begin{array}{ll}
I_2 & 0\\
0 & - I_2
\end{array}\right), 
\quad \alpha_k = \left(
\begin{array}{ll}
0 & \sigma_k\\
 \sigma_k & 0
\end{array}\right), 
$$
for $k = 1,2,3$ and $\sigma_k$ given by
$$\sigma_1 = \left(
\begin{array}{ll}
0 & 1\\
1 & 0
\end{array}\right), 
\, \, \sigma_2 = \left(
\begin{array}{ll}
0 & -i\\
i & 0
\end{array}\right), 
\, \, \sigma_3 = \left(
\begin{array}{ll}
1 & 0\\
0 & -1
\end{array}\right).
$$
Hence, the operator $D_c$ is well defined in $L^2(\mathbb{R}^3, \mathbb{C}^4)$.

Let us consider the nonlinar Dirac equation
\begin{equation}
-i\hbar\frac{\partial\psi}{\partial t} = ic\hbar\sum_{k = 1}^3 \alpha_k \partial_k\psi - mc^2\beta\psi - M(x)\psi + g(x,\psi).
\label{Diracevolucao}
\end{equation}
Assuming that $g(x,e^{i\theta}\psi) = g(x,\psi)$, by the {\it Ansatz} $\psi(t,x) = e^\frac{i\mu  t}{\hbar}u(x)$, one can verify that $\psi(x,t)$ satisfy (\ref{Diracevolucao}) if and only if  $u: \mathbb{R}^3 \to \mathbb{C}^4$ satisfy the following problem
\begin{equation}
-i\hbar\sum_{k = 1}^3 \alpha_k \partial_k u + a\beta u + V(x)u = f(x,u),
\label{Dirac}
\end{equation}
where $a = mc^2$, $V(x) = M(x)/c + \mu I_4$ and $f(x,u) = g(x,u)/c$.

There are many works dedicated to study the Dirac equation (\ref{Dirac}) with the potential $V$ and the nonlinearity $f$ under several different hypotheses. In \cite{Merle}, Merle study the problem (\ref{Dirac}) with a constant potential $V(x) = \omega \in (-a,a)$ and nonlinearity representing the so called Soler model. As far as variational methods are concerned, it seems that Esteban and S\'er\'e in \cite{SereEsteban} were pioneers in using this kind of method to study (\ref{Dirac}).

Motivated by the versatility that variational methods provide, and by the physical appeal of its deduction, some researchers started to work in several generalizations of results which was known to hold to the Schr\"odinger equation, now to the Dirac one. It is important to cite the work of Yanheng Ding and his collaborators, which wrote an extensive list of papers dealing with this subject. The first one seems to be \cite{DingBartsch} where Ding and Bartsch prove the existence and multiplicity of solutions to (\ref{Dirac}), with a scalar potential $V$ and a periodic nonlinearity. Considering matrix potentials, Ding and Ruf in \cite{DingRuf} perform a study about existence and multiplicity of solutions of (\ref{Dirac}) with an asymptotic linear nonlinearity.

In \cite{DingLiu}, Ding and Liu consider the problem (\ref{Dirac}), with $f(x,u) = \nu|u|^{p-2}u$ and prove the existence and concentration, as $\hbar \to 0$, of ground-state semiclassical solutions. In that work they consider the potential $V: \mathbb{R}^3 \to \mathbb{R}$ satisfying the following global condition which was first considered by Rabinowitz in \cite{Rabinowitz},
\begin{description}
\item [$(R)$] $\displaystyle 0 < \inf_{\mathbb{R}^3}V < \liminf_{|x|\to\infty}V$.
\end{description}

In \cite{DingXu}, Ding and Xu improved the result in \cite{DingLiu}, by considering the same problem and proving the same kind of result, but now considering the potential $V$ under a local condition, like the one considered by Del Pino and Felmer in \cite{DelPino}, and nonlinearity in both, the superlinear and the asymptotic linear case. In that work, in a unified approach, the authors performed a penalization technique which resembles \cite{DelPino} in such a way to recover the Cerami compacness condition.

When dealing variationally with (\ref{Dirac}), the most part of the difficulty arises from the spectrum of $D_c = -i\boldsymbol{\alpha} \nabla + a \beta$. In fact, as a consequence of the fact that $\sigma(D_c) = \sigma_c(D_c) = \mathbb{R} \backslash (-a,a)$, the energy functional associated to (\ref{Dirac}) is strongly indefinite, in the sense that his domain contains two infinite dimensional subspaces, where some part of the energy has opposite sign in each of them. By this reason, it becomes difficult to apply the standard minimax theorems like Mountain Pass Theorem, etc., to this kind of functional. Moreover, the standard Nehari manifold are not well defined in this case. 

In this work we study the following version of the Dirac equation
\begin{equation}
-i \boldsymbol{\alpha}\nabla u + a\beta u + V(x)u = K(x)f(|u|)u, \quad \mbox{in $\mathbb{R}^3$,}
\label{P}
\end{equation}
where the potentials $V, K: \mathbb R^3 \rightarrow \mathbb R$ are continuous in $\mathbb R^3$ and are assumed to satisfy the following general conditions. We say that $(V, K) \in \mathcal{K}$ if
\begin{description}
\item[$(VK_0)$] $V(x), K(x) > 0$ for all $x \in \mathbb R^3$, $V, K \in L^{\infty}(\mathbb R^3)$, with $\|V\|_\infty < a$.
\item[$(VK_1)$]   If $(A_n) \subset \mathbb R^3$ is a sequence of Borel sets such that its Lebesgue measure $|A_n| \leq R$, for all $n \in \mathbb N$ and some $R > 0$, then
$$\lim_{r \rightarrow + \infty} \int_{A_n \cap B_r^c(0)} K(x) =0, \quad\hbox{uniformly in $n \in \mathbb N$}.  $$
\end{description}
Furthermore, one of the below conditions occurs
\begin{description}
\item[$(VK_2)$]   $\displaystyle  \frac{K}{V} \in L^{\infty}(\mathbb R^3)$
\end{description}
or
\begin{description}
\item[$(VK_3)$] there exists $q \in (2, 3)$  such that
$$  \frac{K(x)}{V(x)^{3 - q}} \rightarrow 0 \quad\hbox{as $|x| \rightarrow + \infty$}. $$
\end{description}
Moreover, we assume  the following growth conditions at the origin and at infinity for the continuous  function $f: \mathbb R^+ \rightarrow \mathbb R$:
\begin{description}
\item[($f_{1}$)] $\displaystyle \lim_{s \to 0^+}f(s) = 0$;
\item[($f_{2}$)] there exist $c_1, c_2 > 0$ and $p \in (2,3)$ such that $|f(s)s| \leq c_1|s| + c_2|s|^{p-1}$, for all $s \in \mathbb{R}^+$;
\item[($f_{3}$)] $\displaystyle \lim_{t \to \infty}\frac{F(t)}{t^2} = +\infty$;
\item[($f_{4}$)] $f$ is increasing in $\mathbb{R}^+$.
\end{description} 
\begin{rmk}
Note that, from $(f_1)$ and $(f_4)$, it follows that 
$$\frac{1}{2}f(t)t^2 - F(t) \geq 0, \quad \forall t \in\mathbb{R}^+.$$
\label{(f_3)}
\end{rmk}

The main result of this paper is the following theorem.

\begin{theorem}\label{Theorem0}
Suppose that $(V, K) \in \mathcal{K}$ and $f \in C^0(\mathbb{R})$ verifies $(f_1) - (f_4)$. Then,
problem (\ref{P}) possesses a ground state solution, namely, a solution with the lowest energy among all nontrivial ones.
\end{theorem}

This work has been motivated also by the works of Alves and Souto \cite{AlvesSouto} and Figueiredo and Barile \cite{FigueiredoBarile}, where the authors consider semilinear and quasilinear problems, respectively, with the potentials $V$ and $K$ satisfying the same set of assumptions. In order to overcome the apparent lack of compacness, their approach use a compact embedding of the Sobolev space associated to their problems, in some weighted Lebesgue spaces. In fact, since our problem is defined in an unbounded domain, we also have to overcome this difficulty and in order to do so we proceed as \cite{AlvesSouto} and \cite{FigueiredoBarile}.

In spite of the works \cite{AlvesSouto}, \cite{DingLiu}, \cite{DingXu}, where the authors use the Mountain Pass Theorem to get the solution, in this work we follow a different approach that resembles much more that one implemented by Szulkin and Weth in \cite{SzulkinWeth}. In fact we construct the {\it generalized Nehari set} which is going to be a set whose definition resembles the standard Nehari manifold, containing all the nontrivial critical points of the energy functional associated to (\ref{P}). In fact this set has been firstly presented by Pankov in \cite{Pankov} and later thoroughly studied by Szulkin and Weth in \cite{SzulkinWeth}. In \cite{Zhang1,Zhang2}, Zhang, Tang and Zhang succeed in considering the generalized Nehari set in dealing with Dirac equations under different assumptions on the potential and on the nonlinearity. Despite their works, here we do not use the Palais-Smale condition, since in our approach it was enough to prove that the weak limit of a minimizing sequence on the generalized Nehari set, is nontrivial. It is worth pointing out that, in contrast with \cite{Zhang1,Zhang2}, we have not used Ekeland's variational principle since our nonlinearity is just continuous. In fact, in order to use the weak limit of the minimizer sequence of the energy over the generalized Nehari set, we had to apply a Deformation Lemma in an appropiately way.

The paper is organized as follows. In Section 2 we present the variational framework and the compactness result that we use in our approach. In Section 3 we define the generalized Nehari set and prove some of its properties. In the last section we complete the proof of the main result.

\section{Variational framework}

\hspace{0.5cm} Let us consider $D_c = -i \boldsymbol{\alpha} \nabla + a \beta$ which is a self-adjoint operator in $L^2(\mathbb{R}^3,\mathbb{C}^4)$ such that $\mathcal{D}(D_c) = H^1(\mathbb{R}^3,\mathbb{C}^4)$. Since $\sigma(D_c) = \sigma_c(D_c) = (-\infty, a] \cup [a, +\infty)$, it follows that there exists an ortogonal decomposition of $L^2(\mathbb{R}^3,\mathbb{C}^4)$,
$$L^2(\mathbb{R}^3,\mathbb{C}^4) = L^+ \oplus L^-$$
such that $D_c$ is positive definite in $L^+$ and negative definite in $L^-$. Moreover, whenever $u \in L^2(\mathbb{R}^3,\mathbb{C}^4)$, let us denote by $u = u^+ + u^-$ the decomposition of $u$ where $u^+ \in L^+$ and $u^- \in L^-$.

Denoting by $|D_c|$ the absolute value of $D_c$ and by $|D_c|^\frac{1}{2}$ its square root, let $E = \mathcal{D}(|D_c|^\frac{1}{2})$, endowed with the following inner product
$$\langle u,v \rangle := \mbox{Re}\langle |D_c|^\frac{1}{2}u, |D_c|^\frac{1}{2}v \rangle_{L^2}$$
which gives rise to the norm $\|.\| = \langle \cdot, \cdot \rangle^\frac{1}{2}$.

By \cite{LivroDing}[Lemma 7.4], it follows that $E = H^\frac{1}{2}(\mathbb{R}^3,\mathbb{C}^4)$ and $\|\cdot\|$ is equivalent to the usual norm of $H^\frac{1}{2}(\mathbb{R}^3,\mathbb{C}^4)$. Let us remark that  $E$ is continuously embedded into $L^q(\mathbb{R}^3, \mathbb{C}^4)$ for $q \in [2,3)$ (see \cite{Nezza} for example).

Since $E \subset L^2(\mathbb{R}^3,\mathbb{C}^4)$, it follows that
$$E = E^+ \oplus E^-,$$
where $E^+ = E \cap L^+$, $E^- = E \cap L^-$ and the sum is ortogonal with respect to both $\langle \cdot, \cdot \rangle$ and $\langle \cdot, \cdot \rangle_{L^2}$.
\begin{rmk}
Since $\sigma_c(D_c) = (-\infty, a] \cup [a, +\infty)$, it follows that
$$a\|u\|_2^2 \leq \|u\|^2, \quad \forall u \in E.$$
\label{remark0}
\end{rmk}

Associated to problem (\ref{Dirac}) we have the functional $\Phi: E \to \mathbb{R}$ given by
$$\Phi(u) = \frac{1}{2}\left(\|u^+\|^2 - \|u^-\|^2\right) + \frac{1}{2}\int_{\mathbb{R}^3}V(x)|u|^2 dx - \int_{\mathbb{R}^3}K(x)F(|u|)dx.$$
It follows by standard arguments that $\Phi \in C^1(E,\mathbb{R})$. Also, for $u, v \in E$, note that 
\begin{eqnarray*}
\Phi'(u).v & = & \langle u^+, v^+\rangle - \langle u^-, v^-\rangle + \mbox{Re}\int_{\mathbb{R}^3}V(x)u \cdot v dx - \mbox{Re}\int_{\mathbb{R}^3}K(x)f(|u|)u \cdot vdx\\
& = & \mbox{Re}\langle |D_c|^\frac{1}{2}u^+, |D_c|^\frac{1}{2}v^+\rangle_{L^2} - \mbox{Re}\langle |D_c|^\frac{1}{2}u^-, |D_c|^\frac{1}{2}v^-\rangle_{L^2}\\
& & + \mbox{Re}\int_{\mathbb{R}^3}V(x)u \cdot v dx - \mbox{Re}\int_{\mathbb{R}^3}K(x)f(|u|)u \cdot vdx\\
& = & \mbox{Re}\langle u^+, |D_c|v^+\rangle_{L^2} - \mbox{Re}\langle u^-, |D_c|v^-\rangle_{L^2} + \mbox{Re}\int_{\mathbb{R}^3}V(x)u \cdot v dx\\
& &  - \mbox{Re}\int_{\mathbb{R}^3}K(x)f(|u|)u \cdot vdx\\
& = & \mbox{Re}\langle u, D_cv\rangle_{L^2} + \mbox{Re}\int_{\mathbb{R}^3}V(x)u \cdot v dx - \mbox{Re}\int_{\mathbb{R}^3}K(x)f(|u|)u \cdot vdx
\end{eqnarray*}
where $u \cdot v$ denotes the usual inner product in $\mathbb{C}^4$, i.e., $u \cdot v = \sum_{i=1}^4 u_i \overline{v_i}$. In \cite{DingLiu2}[Lemma 2.1] it is proved that critical points of $\Phi$ are weak solutions of (\ref{Dirac}).

To end up this section let us present a compactness result which is going to be used later on.

\begin{prop}
Assume $(V, K) \in \mathcal{K}$. Let $(u_n) \subset E$ be a sequence such that $u_n \rightharpoonup u$ in $E$. Then, 
\begin{description}
\item [$i)$] if $(VK_2)$ holds,  then $\displaystyle \int_{\mathbb{R}^3}K(x)|u_n|^q dx \to \int_{\mathbb{R}^3}K(x)|u|^q dx$, for all $q \in
(2,3)$;
\item [$ii)$] if $(VK_3)$ holds,  then $\displaystyle \int_{\mathbb{R}^3}K(x)|u_n|^q dx \to \int_{\mathbb{R}^3}K(x)|u|^q dx$.
\end{description}
\label{embedding}
\end{prop}

\begin{dem}  In order to prove the first item, assume that $(VK_2)$ holds. Fixed $q \in
(2, 3)$ and $\epsilon > 0$, there exist $0 < t_0 <
t_1$ and a positive constant $C > 0$ such that
\begin{equation*}
K(x)|t|^q \leq \epsilon\,  C 	\,  (V(x) |t|^2 + |t|^3) + C \, 
K(x) \, \chi_{[t_0, t_1]}(|t|) |t|^3, \quad\hbox{ for all $ t \in
\mathbb R$.}
\end{equation*}
Then, denoting $Q(u)= \displaystyle\int_{\mathbb R^3} V(x) |u|^2
dx + \displaystyle\int_{\mathbb R^3} |u|^3 dx $ and
$A= \left\{ x \in \mathbb R^3: t_0 \leq |u(x)| \leq t_1   \right\}$,
we have that
\begin{equation}\label{2.5}
\int_{B_r^c(0)} K(x) |u|^q dx \leq \epsilon \,  C Q(u) + C
\int_{A \cap B_r^c(0)} K(x) dx, \quad\hbox{for all $u \in E$.}
\end{equation}
Since $(u_n)$ is a weakly convergent sequence, by Banach-Steinhaus Theorem it is bounded in $E$. By the continuous embedding $E \hookrightarrow L^2(\mathbb{R}^3,\mathbb{C}^4)$, $E \hookrightarrow L^3(\mathbb{R}^3,\mathbb{C}^4)$ and the fact that $V \in L^\infty(\mathbb{R}^3)$, there exists $C' > 0$ such that
$$ \int_{\mathbb R^3} V(x) |u_n|^2 dx \leq C' \quad \hbox{and} \quad \int_{\mathbb R^3} |u_n|^3 dx \leq C', \quad\hbox{for all $n \in \mathbb N$}.$$
Then there exists $C'' > 0$ such that $Q(u_n) \leq C''$ for all $n \in \mathbb{N}$. On the other hand, denoting $A_n=
\{ x \in \mathbb R^3: t_0 \leq |u_n(x)| \leq t_1 \}$, it follows
that
$$ t_0^3 |A_n| \leq \int_{A_n} |u_n|^3 dx \leq C', \quad\hbox{for any $n \in \mathbb N$},$$
and then $\sup_{n \in \mathbb N} |A_n| < + \infty$. Consequently, from $(VK_1)$ there exists a positive radius $r > 0$ large enough such that
\begin{equation}\label{2.6}
\int_{A_n \cap B_r^c(0)} K(x) dx < \frac{\epsilon}{t_1^3} \quad\hbox{for all $n \in \mathbb N$}.
\end{equation}
By \eqref{2.5} and \eqref{2.6} it follows that
\begin{eqnarray}\label{2.7}
\int_{B_r^c(0)} K(x) |u_n|^q dx &\leq&CC''\epsilon 
+ C \int_{A_n \cap B_r^c(0)} K(x) dx \nonumber \\
&\leq &
(CC'' + C/{t_1^3}) \epsilon \quad\hbox{for all $n \in \mathbb N$}.
\end{eqnarray}
Since $q \in (2, 3)$ and $K$ is a continuous
function, from Sobolev embeddings we have that
\begin{equation}\label{2.8}
\lim_{n \rightarrow + \infty } \int_{B_r(0)} K(x) |u_n|^q dx
= \int_{B_r(0)} K(x) |u|^q dx.
\end{equation}
Then, from \eqref{2.7} for $\epsilon > 0$ small enough and \eqref{2.8} it follows that
\begin{equation*}
\lim_{n \rightarrow + \infty } \int_{\mathbb R^3} K(x) |u_n|^q
dx = \int_{\mathbb R^3} K(x) |u|^q dx
\end{equation*}
which concludes the proof in this case.

\medskip

Now suppose that $(VK_3)$ holds. Let us define for each $x \in \mathbb R^3$ fixed, the function
$$g(t)= V(x) t^{2-q} + t^{3-q}, \quad\hbox{for every $t > 0$}. $$
Since its minimum value is $C_q V(x)^{3-q}$ where 
$$C_q = \frac{1}{3-q}\left(\frac{q-2}{3-q}\right)^{2-q},$$
we have that
$$C_q V(x)^{3-q} \leq V(x) t^{2-q} + t^{3-q}, \quad\hbox{for every $x \in \mathbb R^3$ and $t > 0$}. $$
Combining this inequality  with $(VK_3)$, for any $\epsilon > 0$ 
there exists a positive radius $r > 0$ sufficiently large
such that
\begin{equation*}
 K(x) |t|^q \leq \epsilon \, C'_q (V(x) |t|^2 + |t|^3), \quad\hbox{for every $t \in \mathbb R$ and $|x| > r$,}
 \end{equation*}
where $C'_q = C_q^{-1}$,  from which it follows that
$$ \displaystyle\int_{B_r^c(0)} K(x) |u|^q dx \leq \epsilon C'_q \displaystyle\int_{B_r^c(0)} (V(x) |u|^2 + |u|^3) dx , \quad\hbox{for all $u \in E$}. $$

Then, for the sequence $(u_n)$ of the statement, again by Banach Steinhaus Theorem, Sobolev embeddings and the boundedness of $V$, there exists $C' > 0$ such that
$$ \int_{\mathbb R^3} V(x) |u_{n}|^2 \leq C' \quad \hbox{and} \quad \int_{\mathbb R^3} |u_{n}|^3 dx \leq C', \quad\hbox{for all $n \in \mathbb N$}, $$
and then
\begin{eqnarray}\label{2.10}
\int_{B_r^c(0)} K(x) |u_n|^q dx \leq 2C' C'_q \epsilon \quad \forall n \in \mathbb N.
\end{eqnarray}
Since $q \in (2, 3)$ and $K$ is a continuous function,
from Sobolev embeddings on bounded domains, we have that
\begin{equation}\label{2.11}
\lim_{n \rightarrow + \infty } \int_{B_r(0)} K(x) |u_n|^q dx = \int_{B_r(0)} K(x) |u|^q dx.
\end{equation}
Then, from \eqref{2.10} for $\epsilon > 0$ small enough  and \eqref{2.11}, it holds that
\begin{equation*}
\lim_{n \rightarrow + \infty } \int_{\mathbb R^3} K(x) |u_n|^qdx = \int_{\mathbb R^3} K(x) |u|^qdx.
\end{equation*}
which concludes the proof.
\end{dem}

\begin{rmk}
It follows from the last proposition, $(f_2)$ and from Lebesgue Dominated Convergence Theorem, that, if $u_n \rightharpoonup u$ in $E$, then
$$
\int_{\mathbb{R}^3}K(x)F(|u_n|)dx \to \int_{\mathbb{R}^3}K(x)F(|u|)dx, \quad \mbox{as $n \to \infty$.}
$$
\label{remark}
\end{rmk}

\section{Generalized Nehari set and its properties} 

\hspace{0.5cm} Let us consider the following set
$$\mathcal{M} = \{u \in E \backslash E^-; \, \Phi'(u)u = 0 \, \, \mbox{and} \, \, \Phi'(u)v = 0, \, \mbox{for all $v \in E^-$}\},$$
which has been introduced by Pankov in \cite{Pankov}, deeply studied by Szulkin and Weth in \cite{SzulkinWeth} and been called generalized Nehari set. 

\begin{rmk}
Note that, if $u \in E$, $u \neq 0$ and $\Phi'(u) = 0$, then by Remark \ref{(f_3)},
$$\Phi(u) = \Phi(u) - \frac{1}{2}\Phi'(u)u = \int_{\mathbb{R}^3}K(x)\left(\frac{1}{2}f(|u|)|u|^2 - F(|u|)\right)dx > 0.$$
Note also that, for $u \in E^-\backslash\{0\}$, it follows by $(VK_0)$ and Remark \ref{remark0} that
\begin{eqnarray*}
\Phi(u) & = & -\frac{1}{2}\|u\|^2 + \frac{1}{2}\int_{\mathbb{R}^3}V(x)|u|^2dx - \int_{\mathbb{R}^3}K(x)F(|u|)dx\\
& \leq & -\frac{a}{2}\|u\|_2^2 + \frac{\|V\|_\infty}{2}\|u\|_2^2\\
& < & 0.
\end{eqnarray*}
Then all nontrivial critical points of $\Phi$ belong to $E \backslash E^-$ and then $\mathcal{M}$ contains all nontrivial critical points. Then, if $u \in \mathcal{M}$ is a critical point of $\Phi$, then it is going to have the lowest energy among all nontrivial critical points, justifying we calling it a ground state solution of (\ref{Dirac}).
\end{rmk}

Let us follow the notation stablished by Szulkin and Weth in \cite{SzulkinWeth} and, for $u \in E\backslash E^-$, let us denote $\hat{E}(u) = \mathbb{R}^+u^+ \oplus E^-$. Also, for each $u \in E\backslash E^-$, let us define $\gamma_u: \mathbb{R}^+ \times E^- \to \mathbb{R}$ given by
$$\gamma_u(t,v) = \Phi(tu^+ + v),$$
for all $(t,v) \in \mathbb{R}^+ \times E^-$ and note that $\gamma_u \in C^1(\mathbb{R}^+ \times E^-, \mathbb{R})$.
Note also that
\begin{equation}
\frac{\partial}{\partial t}\gamma_u(t,v) = \Phi'(tu^+ + v)u^+
\label{Nehari1}
\end{equation}
and
\begin{equation}
\frac{\partial}{\partial v}\gamma_u(t,v)w = \Phi'(tu^+ + v)w, \quad \forall w \in E^-.
\label{Nehari2}
\end{equation}

\begin{lem}
$(t,v) \in \mathbb{R}^+ \times E^-$ is a critical point of $\gamma_u$ if and only if $t u^+ + v \in \mathcal{M}$.
\label{lemmaNehari1}
\end{lem}
\begin{dem}
In fact, suppose that 
$$\frac{\partial}{\partial t}\gamma_u(t,v) = 0$$
and
$$\frac{\partial}{\partial v}\gamma_u(t,v) = 0, \quad \mbox{in ${(E^-)}'$.}$$
By (\ref{Nehari1}) and (\ref{Nehari2}), we have that
$$\Phi'(t u^+ + v)(t u^+ + v) = t\Phi'(t u^+ + v)u^+ + \Phi'(t u^+ + v)v = 0.$$
Moreover, again by (\ref{Nehari2}),
$$\Phi'(tu^+ + v)w = 0, \quad \forall w \in E^-,$$
and then $tu^+ + v \in \mathcal{M}$.

Conversely, if $tu^+ + v \in \mathcal{M}$, where $(t,v) \in \mathbb{R}^+ \times E^-$, then
$$\frac{\partial}{\partial v}\gamma_u(t,v)w = \Phi'(tu^+ + v)w = 0, \quad \forall w \in E^-$$
and also, taking into account the last expression,
$$\Phi'(tu^+ + v)(tu^+ + v) = 0 \Leftrightarrow t\Phi'(tu^+ + v)u^+ = 0$$
which implies that
$$
\frac{\partial}{\partial t}\gamma_u(t,v) = 0.
$$
Then $(t,v)$ is a critical point for $\gamma_u$.
\end{dem}

\begin{lem}
For each $u \in E\backslash E^-$, there exists $t_u u^+ + v_u \in \hat{E}(u)$, such that
$$\Phi(t_u u^+ + v_u) = \max_{t \geq 0,\, v \in E^-} \Phi(t u^+ + v).$$
Moreover, $t_u u^+ + v_u \in \mathcal{M}$.
\label{lemma2.2}
\end{lem}
\begin{dem}
Note that by Lemma \ref{lemmaNehari1}, if the maximum exists then it is going to belong to $\mathcal{M}$.

Since $\hat{E}(u) = \hat{E}\left(u^+/\|u^+\|\right)$ we can assume without lack of generality that $u \in E^+$ and $\|u^+\| = 1$. 
Note that by $(f_1)$, $\Phi(tu) > 0$ for all $t > 0$ sufficiently small.

Now let us prove that there exists $R > 0$ such that 
\begin{equation}
\Phi(w) \leq 0 \quad \forall w \in \hat{E}(u)\backslash B_R(0).
\label{limitacao}
\end{equation}

Suppose by contradiction that there exist $(w_n) \subset \hat{E}(u)$ such that $\|w_n\| \to +\infty$ and $\Phi(w_n) > 0$. By doing $w_n = t_n u + v_n$ where $t_n \geq 0$ and $v_n \in E^-$, let us define
$$
\overline{w}_n = \frac{w_n}{\|w_n\|} = \frac{t_n}{\|w_n\|}u + \frac{v_n}{\|w_n\|} =: s_n u + \overline{v_n}.
$$
Note that
\begin{equation}
1 = \|\overline{w_n}\|^2 = s_n^2 + \|\overline{v}_n\|^2.
\label{Nehari3}
\end{equation}
Also, since $F(t) \geq 0$ in $\mathbb{R}^+$, by Remark \ref{remark0}
\begin{eqnarray*}
0 & < & \frac{\Phi(w_n)}{\|w_n\|^2}\\
& = & \frac{1}{2}(s_n^2 - \|\overline{v}_n\|^2) + \frac{1}{2}\int_{\mathbb{R}^3}V(x)|s_nu + \overline{v}_n|^2dx - \int_{\mathbb{R}^3}\frac{F(|w_n|)}{\|w_n\|^2}dx\\
& \leq & \frac{1}{2}(s_n^2 - \|\overline{v}_n\|^2) + \frac{1}{2}\|V\|_\infty\|s_nu + \overline{v}_n\|_2^2\\
& = & \frac{1}{2}(s_n^2 - \|\overline{v}_n\|^2) + \frac{1}{2}\|V\|_\infty(s_n^2\|u\|_2^2 + \|\overline{v}_n\|_2^2)\\
& \leq & \frac{1}{2}\left[s_n^2\left(1+ \frac{\|V\|_\infty}{a}\right) - \|\overline{v}_n\|^2\left(1- \frac{\|V\|_\infty}{a}\right)\right]
\end{eqnarray*}
which together with (\ref{Nehari3}) impliy that
\begin{equation}
\left(\frac{a - \|V\|_\infty}{a + \|V\|_\infty}\right)  \|\overline{v}_n\|^2 \leq s_n^2 = 1- \|\overline{v}_n\|^2.
\label{Nehari4}
\end{equation}
Then, by (\ref{Nehari4}), it follows that
$$
0 \leq \|\overline{v}_n\|^2 \leq \frac{a + \|V\|_\infty}{2a} \quad \mbox{and} \quad 0 < \frac{a-\|V\|_\infty}{2a} \leq s_n \leq 1, \quad \forall n \in \mathbb{N}.
$$
Then, we can assume that, up to a subsequence $\overline{v}_n \rightharpoonup \overline{v}$ and $s_n \to s_0 \neq 0$. Hence
$$
\overline{w_n} \rightharpoonup \overline{w} = s_0 u + \overline{v} \neq 0.
$$

Let us denote by $\Gamma = \{x \in \mathbb{R}^3; \, \overline{w}(x) \neq 0\}$ and note that $|\Gamma| > 0$.
On one hand we have that
$$0 < \frac{\Phi(w_n)}{\|w_n\|^2}, \quad \forall n \in \mathbb{N},$$
while in the other, by $(f_3)$ and Remark \ref{remark0}
\begin{eqnarray*}
\limsup_{n \to \infty} \frac{\Phi(w_n)}{\|w_n\|^2} & = & \frac{s_0^2}{2} - \frac{1}{2}\liminf_{n \to \infty}\|\overline{v}_n\|^2 + \frac{1}{2}\liminf_{n \to \infty}\int_{\mathbb{R}^3}V(x)\frac{|w_n|^2}{\|w_n\|^2}dx\\
& &  - \liminf_{n \to \infty}\int_{\mathbb{R}^3}K(x)\frac{F(|w_n|)}{\|w_n\|^2}dx\\
& \leq & \frac{1}{2}\left(s_0^2 - \|\overline{v}\|^2\right) + \frac{\|V\|_\infty}{a} - \liminf_{n \to \infty}\int_{\Gamma}K(x)\frac{F(|w_n|)}{\|w_n\|^2}dx\\
& = & - \infty,
\end{eqnarray*}
which give us a clear contradiction.

Since $\Phi$ is bounded from above in $\hat{E}(u)$, let us take a maximizing sequence $(u_n) \subset \hat{E}(u)$ such that
$$
\lim_{n \to \infty}\Phi(u_n) = \beta := \max_{\hat{E}(u)}\Phi.
$$
Since $0 < \beta < +\infty$ and by (\ref{limitacao}), it follows that there exists $R > 0$ such that $\|u_n\| \leq R$, for all $n \in \mathbb{N}$. Writing $u_n = r_n u + u_n^-$ note that
$$
R^2 \geq \|u_n\|^2 = r_n^2 + \|u_n^-\|^2,
$$
which implies that both $(r_n) \subset \mathbb{R}$ and $(u_n^-) \subset E^-$ are bounded. Then, up to a subsequence, $r_n \to r_0$ and $u_n^- \rightharpoonup u_0^-$. Then it follows that
$$
u_n = r_n u + u_n^- \rightharpoonup r_0 u + u_0^- =: u_0 \in \hat{E}(u).
$$
Then, the last informations, since $F \geq 0$, Fatou Lemma imply that
\begin{eqnarray*}
- \Phi(u_0) & = & \frac{1}{2}(\|u_0^-\|^2 - r_0^2) - \frac{1}{2}\int_{\mathbb{R}^3}V(x)|u_0|^2 dx+\int_{\mathbb{R}^3}K(x)F(|u_0|)dx\\
& \leq & \liminf_{n \to \infty}\left[\frac{1}{2}(\|u_n^-\|^2 - r_0^2) + \int_{\mathbb{R}^3}K(x)F(|u_n|)dx\right]\\
& = & \liminf_{n \to \infty} -\Phi(u_n) = - \beta.
\end{eqnarray*}
Then $\Phi(u_0) \geq \beta$ and hence $\Phi(u_0) = \beta$, with $u_0 \in \hat{E}(u)$.
\end{dem}

The next result is a technical lemma that is going to be used to prove the uniqueness of the maximum point of $\Phi|_{\hat{E}(u)}$. Its proof follows the arguments as in \cite{SzulkinWeth}, but we present it here for the sake of completeness.

\begin{lem}
Let $t \in \mathbb{R}$, $t \geq 0$ and $u,v \in \mathbb{C}^4$ such that $v \neq 0$. Then
$$
\mbox{Re}f(|u|)u \cdot \left(\frac{t^2}{2}u - \frac{u}{2} + tv \right) + F(|u|) - F(|tu+v|) < 0,
$$
where $u \cdot v$ denotes the usual inner product in $\mathbb{C}^4$.
\label{lemma3.3}
\end{lem}
\begin{dem}
Let $t, u$ and $v$ as in the statement. Define
$$h(t):= \mbox{Re} f(|u|)u \cdot \left(\frac{t^2}{2}u - \frac{u}{2} + tv \right) + F(|u|) - F(|tu+v|)$$
and note that we have to prove that $h(t) < 0$.

Note that if $u = 0$, then $h(t) = -F(|v|) < 0$ by Remark \ref{(f_3)}. Then we can assume $u \neq 0$.
We have to consider two cases. Suppose first that $\mbox{Re}\left( u \cdot (tu+v)\right) \leq 0$. By using Remark \ref{(f_3)},
\begin{eqnarray*}
h(t) & < & \mbox{Re} f(|u|)u \cdot \left(\frac{t^2}{2}u - \frac{u}{2}\right)u + \frac{1}{2}f(|u|)u\cdot u + \mbox{Re}\, \,  tf(|u|)u\cdot v - F(|tu+v|)\\
& = &  \frac{t^2}{2}f(|u|)u\cdot u + \mbox{Re} \, \, tf(|u|)u\cdot v - F(|tu + v|)\\
& = & -\frac{t^2}{2}f(|u|)u\cdot u + \mbox{Re}\, \, tf(|u|)u\cdot(tu + v) - F(|tu + v|) \leq 0.
\end{eqnarray*}
Now suppose that $\mbox{Re} \left(u\cdot (tu+v)\right) > 0$, and note that
$$h(0) = -\frac{1}{2}f(|u|)u\cdot u + F(|u|) < 0$$
by Remark \ref{(f_3)}. Moreover note that $\lim_{t \to +\infty}h(t) = - \infty$.
Note also that
$$h'(t) = \mbox{Re}\left(u\cdot (tu+v)\right) \left(f(|u|) - f(|tu+v|)\right).$$
Assuming that there exists a maximum point $t_0 \geq 0$ such that $h(t_0) \geq 0$, then, since $h'(t_0) = 0$ and $\mbox{Re} \left(u\cdot (tu+v)\right) > 0$, it follows by $(f_4)$ that $|u| = |t_0u + v|$. 
Then, observing that $u\cdot u = \displaystyle \max_{v \in \mathbb{C}^4}\mbox{Re} \, \, u \cdot v$, we have that
\begin{eqnarray*}
h(t_0) & = & \mbox{Re}\left(f(|u|)u \cdot \left(\frac{t_0^2}{2}u - \frac{u}{2}\right) + t_0f(|u|)u\cdot v\right)\\
& = & \mbox{Re}\left(-\frac{t_0^2}{2}f(|u|)u^2 - \frac{1}{2}f(|u|)u^2 + t_0f(|u|)u(t_0u+v)\right)\\
& < & -\frac{t_0^2}{2}f(|u|)u^2 - \frac{1}{2}f(|u|)u^2 + t_0f(|u|)u^2\\
& = & -\frac{(t_0-1)^2}{2}f(|u|)u^2\\
&  \leq & 0,
\end{eqnarray*}
which leads us to a contradiction. 

Therefore, in any case $h(t) < 0$ for all $t \in \mathbb{R}^+$.
\end{dem}

\begin{lem}
For each $u \in \mathcal{M}$, we have that
$$\Phi(tu + v) < \Phi(u) \quad \forall (tu+v) \in \hat{E}(u).$$
\end{lem}
\begin{dem}
Let $u \in \mathcal{M}$, $t \geq 0$ and $v \in E^-$. Note that, since $\Phi'(u)w = 0$ for all $w \in \hat{E}(u)$, we have that
\begin{eqnarray*}
\Phi(tu + v) - \Phi(u) & = & \frac{1}{2}\left(\langle D_c( tu + v), tu + v\rangle_{L^2} - \langle D_cu, u \rangle_{L^2}\right)\\
& & + \frac{1}{2}\int_{\mathbb{R}^3}V(x)\left(|tu + v|^2 - |u|^2\right)dx\\
& & + \int_{\mathbb{R}^3}K(x)\left(F(|u|) - F(|tu+v|)\right)dx\\
& = & \mbox{Re}\left<D_cu, \frac{t^2u}{2} - \frac{u}{2} + tv\right>_{L^2} - \frac{1}{2}\|v\|^2\\
& & +  \mbox{Re}\int_{\mathbb{R}^3}V(x)u\cdot\left(\frac{t^2u}{2} - \frac{u}{2} + tv\right)dx - \frac{1}{2}\int_{\mathbb{R}^3}V(x)|v|^2dx\\
& & + \int_{\mathbb{R}^3}K(x)\left(F(|u|) - F(|tu+v|)\right)dx\\
& = &  \int_{\mathbb{R}^3}\mbox{Re}\left[K(x)f(|u|)u\cdot\left(\frac{t^2}{2}u - \frac{u}{2} + tv \right) + F(|u|) - F(|tu+v|)\right]dx\\
& & - \frac{1}{2}\|v\|^2 - \frac{1}{2}\int_{\mathbb{R}^3}V(x)|v|^2\\
& < & 0,
\end{eqnarray*}
by Lemma \ref{lemma3.3}, since $\displaystyle \frac{t^2u}{2} - \frac{u}{2} + tv \in \hat{E}(u)$. Then the result follows.
\end{dem}

We can summarize the last three results in the following proposition.
\begin{prop}
For each $u \in E\backslash E^-$, there exists a unique $t_u > 0$ and $v_u \in E^-$ such that $t_u u^+ + v_u \in \mathcal{M}$. Moreover, if $u \in \mathcal{M}$, then $t_u = 1$ and $v_u = u^-$.
\label{prop3.1}
\end{prop} 

Now, a picture is in order to clarify the main properties of the generalized Nehari set.
\begin{figure}[!htb]
\center
\includegraphics{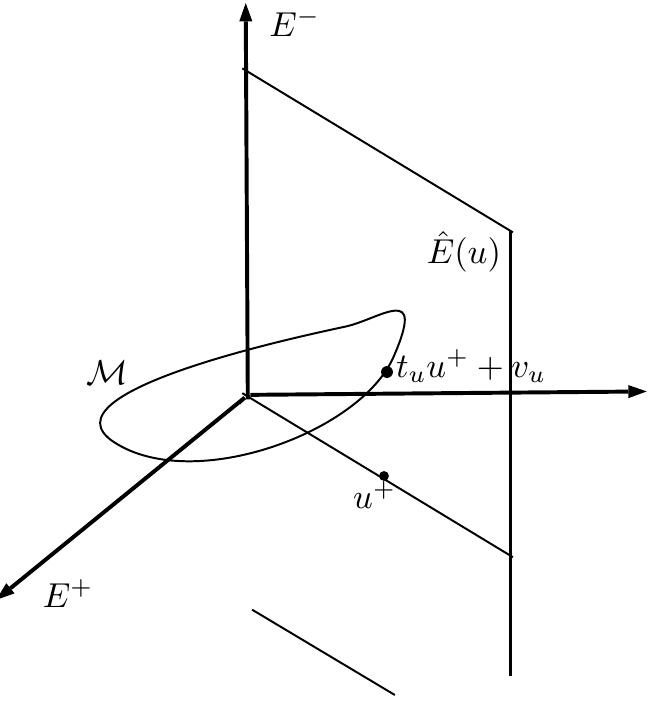}
\caption{Generalized Nehari set}
\end{figure}

Clearly the last result implies that there exists a bijective map between $\mathcal{M}$ and $\{u \in E^+; \|u\| = 1\}$. Also, it holds the following result.

\begin{lem}
There exists $\delta > 0$ such that, for all bounded subset $S \subset \mathcal{M}$, we have that
$$
\|u^+\| \geq \delta, \quad \forall u \in S.
$$
\end{lem}
\label{cotainferiorNehari}
\begin{dem}
First of all, let us prove that there exists $\alpha > 0$ such that $\Phi(u) \geq \alpha$, for all $u \in \mathcal{M}$. Let us remark before it that there exist $\rho, \alpha > 0$ such that
$$\Phi(u^+) \geq \alpha \quad \forall u^+ \in E^+ \cap B_\rho(0).$$
In fact, just note that by $(f_1)$ and $(f_2)$, for all $\epsilon > 0$, there exists $A_\epsilon > 0$ such that
$$|f(s)s| \leq \epsilon |s| + A_\epsilon |s|^{p-1}, \quad \forall s \in \mathbb{R}^+,$$
where $p \in (2,3)$ is like in $(f_2)$.
Then using this inequality and Sobolev embeddings we have that
\begin{eqnarray*}
\Phi(u^+) & = & \frac{1}{2}\|u^+\|^2 + \frac{1}{2}\int_{\mathbb{R}^3}V(x)|u^+|^2dx - \int_{\mathbb{R}^3}K(x)F(|u^+|)dx\\
& \geq & \left(\frac{1}{2} - C\epsilon\right)\|u^+\|^2 - CA_\epsilon \|u^+\|^p\\
& = & \|u^+\|^2\left(\frac{1}{2} - C\epsilon - CA_\epsilon\|u^+\|^{p-2}\right)\\
& \geq & \alpha,
\end{eqnarray*}
just by choosing $0 < \epsilon < \frac{1}{2C}$, for all $u^+ \in E^+$ such that $\|u^+\| = \rho$, where $\rho = \left(\frac{1}{CA_\epsilon}\left(\frac{1}{2} - C\epsilon\right)\right)^{\frac{1}{p-2}}$ and $\alpha = \rho^2\left(\frac{1}{2} - C\epsilon - CA_\epsilon\rho^{p-2}\right)$.

Now just note that by Lemma \ref{lemma2.2}, for all $u \in \mathcal{M}$
$$\Phi(u) \geq \Phi\left(\rho\frac{u^+}{\|u^+\|}\right) \geq \alpha.$$

In order to finish the proof, let us assume by contradiction that there exists a bounded sequence $(u_n) \subset \mathcal{M}$, such that $u_n^+ \to 0$ in $E$. Since $u_n \in \mathcal{M}$, it follows that
\begin{equation}
\|u_n^+\|^2 - \|u_n^-\|^2 + \int_{\mathbb{R}^3}V(x)u_n \cdot u_ndx = \int_{\mathbb{R}^3}K(x)f(|u_n|)u_n \cdot u_ndx,
\label{eqNehari1}
\end{equation}
and
\begin{equation}
- \|u_n^-\|^2 + \mbox{Re}\int_{\mathbb{R}^3}V(x)u_n \cdot u_n^-dx = \mbox{Re}\int_{\mathbb{R}^3}K(x)f(|u_n|)u_n \cdot u_n^-dx.
\label{eqNehari2}
\end{equation}
Then, by using (\ref{eqNehari2}) and since $u_n^+ \to 0$, by Sobolev embeddings and H\"older inequality we have that
\begin{eqnarray*}
\alpha & \leq & \Phi(u_n)\\
& = & \frac{1}{2}\|u_n^+\|^2 + \frac{1}{2}\mbox{Re}\int_{\mathbb{R}^3}V(x)u_n \cdot u_n^+dx\\
& & + \frac{1}{2}\mbox{Re}\int_{\mathbb{R}^3}K(x)f(|u_n|)u_n \cdot u_n^-dx - \int_{\mathbb{R}^3}K(x)F(|u_n|)dx\\
& \leq & o_n(1) +  \frac{1}{2}\int_{\mathbb{R}^3}K(x)f(|u_n|)|u_n^-|^2dx.
\end{eqnarray*}
Above we have used the fact that by $(f_1)$ and $(f_2)$, for all $\epsilon > 0$, there exists $A_\epsilon > 0$ such that
$$|f(s)| \leq \epsilon + A_\epsilon |s|^{p-2}, \quad \forall s \in \mathbb{R}^+,$$
together with H\"older inequality.
Then, again by (\ref{eqNehari2}), since $u_n^+ \to 0$ in $E$, it follows that
\begin{eqnarray*}
2\alpha + o_n(1)& \leq & \int_{\mathbb{R}^3}K(x)f(|u_n|)|u_n^-|^2dx\\
& = & - \|u_n^-\|^2 + \mbox{Re}\int_{\mathbb{R}^3}V(x)u_n \cdot u_n^-dx\\
& = & - \|u_n^-\|^2 + \int_{\mathbb{R}^3}V(x)|u_n^-|^2dx + o_n(1)\\
& \leq & \|V\|_\infty\|u_n^-\|^2_2 - \|u_n^-\|^2 + o_n(1)\\
& \leq & \left(\|V\|_\infty - a\right)\|u_n^-\|^2_2 + o_n(1)\\
& < & 0,
\end{eqnarray*}
where we have used Remark \ref{remark0} and $(VK_0)$ in the last inequalities. This contradiction prove the result.
\end{dem}

Note that $\Phi$ is bounded from bellow in $\mathcal{M}$, since if $u \in \mathcal{M}$, then $\Phi(u) \geq \Phi(0) = 0$ (since $0 \in \hat{E}(u)$). 

To end up this section, let us prove that if the minimum of $\Phi$ on $\mathcal{M}$ is achieved in some $u \in \mathcal{M}$, then in fact $u$ is a critical point of $\Phi$. This follow from a Deformation Lemma (see \cite{Willem}) and is going to be proved in the next result.

\begin{prop}
If $u_0 \in \mathcal{M}$ is such that
$$
\Phi(u_0) = \min_{\mathcal{M}}\Phi,
$$
then $\Phi'(u_0) = 0$.
\end{prop}
\begin{dem}

Suppose by contradiction that $\Phi'(u_0) \neq 0$. By the continuity of $\Phi'$, it follows that there exist $\kappa, \lambda > 0$ such that

\begin{equation}
\|\Phi'(u)\|_{E^*} \geq \lambda, \quad \forall u \in B_\kappa(u_0).
\label{betaestimate}
\end{equation}

\begin{figure}[!htb]
\center
\includegraphics{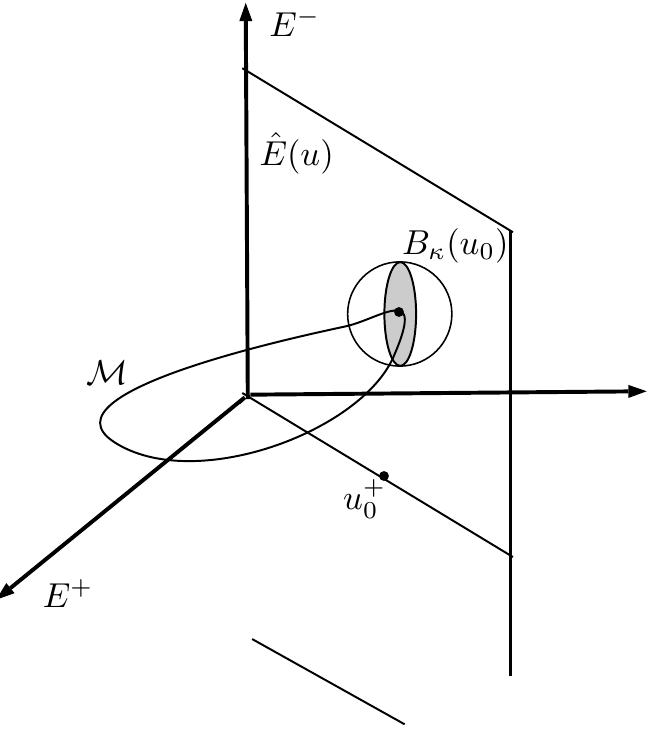}
\end{figure}

Let us define $T: \hat{E}(u_0) \to \mathbb{R}^+ \times E^-$ given by $T(tu_0^+ + v) = (t,v)$. Denoting $B_{\kappa/4} = B_{\kappa/4}(u_0)\cap \hat{E}(u_0)$, let us define $g = T^{-1}|_{T(B_{\kappa/4})}$, i.e., $g : T(B_{\kappa/4}) \to \mathbb{R}^+ \times E^-$ is given by
$$g(t,v) = tu_0^+ + v, \quad \mbox{for $(t,v) \in T(B_{\kappa/4})$.}$$

Denoting $\displaystyle c := \Phi(u_0) = \inf_{\mathcal{M}}\Phi$, note that by Proposition \ref{prop3.1},
$$
\Phi(g(t,v)) < c, \quad \forall (t,v) \neq (1,u^-).
$$
Moreover, since $T$ is an homeomorphism, note that $\partial T(B_{\kappa/4}) = T(\partial B_{\kappa/4})$ and also
$$
\max_{(t,v) \in \partial T(B_{\kappa/4})}\Phi(g(t,v)) = c_0 < c.
$$

Now let us use the Deformation Lemma (see \cite{Willem}[Lemma 2.3]), with 
$$
0 < \epsilon < \min\left\{\frac{c-c_0}{2},\frac{\lambda\kappa}{32}\right\}
$$
and $\delta = \kappa/4$. Then it follows that there exists an homeomorphism $\eta: E \to E$ such that
\begin{itemize}
\item [$i)$] $\eta(x) = x$ for all $x \not \in  \Phi^{-1}([c - 2\epsilon,c + 2\epsilon]) \cap B_\kappa(u_0)$;
\item [$ii)$] $\eta(\Phi_{c+\epsilon}\cap B_{\kappa/2}(u_0)) \subset \Phi_{c-\epsilon}$;
\item [$iii)$] $\Phi(\eta(x)) \leq \Phi(x)$, for all $x \in E$.
\end{itemize}

Let us define now $h: T(B_{\kappa/4}) \to E$ by $h(t,v) = \eta(g(t,v))$ and, for each $w \in E^-$, two functions, $\Psi^w_0, \Psi^w_1 : T(B_{\kappa/4}) \to \mathbb{R}$ by
$$\Psi^w_0(t,v) = (\Phi'(g(t,v))g(t,v), \Phi'(g(t,v))w)$$
and
$$\Psi^w_1(t) = (\Phi'(h(t,v))h(t,v), \Phi'(h(t,v))w).$$

Since for $\displaystyle (t,v) \in T(\partial B_{\kappa/4})$, $\Phi(g(t,v)) \leq c_0 < c - 2\epsilon$, then 
$$h(t,v) = \eta(g(t,v)) = g(t,v), \quad \mbox{for $(t,v) \in T(\partial B_{\kappa/4})$.}$$
Hence
\begin{equation}
\Psi^w_0(t,v) = \Psi^w_1(t,v), \quad \forall (t,v) \in T(\partial B_\kappa).
\label{fronteira}
\end{equation}

By Degree Theory, $d(\Phi^w_0,T(B_{\kappa/4}),(0,0)) = 1$ and, taking into account (\ref{fronteira}), we have that $d(\Phi^w_1,T(B_{\kappa/4}),(0,0)) = 1$. Then we have that there exists $(t_0,v_0) \in T(B_{\kappa/4})$ such that $h(t_0,v_0) \in \mathcal{M}$.
This implies that 
$$c \leq \Phi(h(t_0,v_0)) = \Phi(\eta(g(t_0,v_0))).$$
But note that $\Phi(g(t_0,v_0)) < c + \epsilon$ and also $g(T(B_{\kappa/4})) = B_{\kappa/4} \subset B_{\kappa/2}(u_0)$. Then, by $ii)$
$$ \Phi(\eta(g(t_0,v_0))) < c - \epsilon$$
which contradicts the last inequality. Then the result follows.
\end{dem}

\section{Proof of the main result}

\hspace{0.5cm} Let us take a minimizing sequence $(u_n) \in \mathcal{M}$ such that 
$$\Phi(u_n) \to \min_{w \in \mathcal{M}} \Phi(w) =: c.$$ 

It is an easy matter to prove that $(u_n)$ is a bounded sequence. Then there exists $u \in E$ such that $u_n \rightharpoonup u$ in $E$ up to a subsequence. Note also that $u_n^+ \rightharpoonup u^+$ and $u_n^- \rightharpoonup u^-$ by the continuity of the projections maps. 

Now let us prove that $u^+ \neq 0$. On the contrary,  $u_n^+ \rightharpoonup 0$. Then, for all $t \geq 0$, since $\Phi(u_n) = \max_{\hat{E}(u_n)}\Phi$ and by Lemma \ref{cotainferiorNehari} and Remark \ref{remark}, it follows that
\begin{eqnarray*}
c + o_n(1) & \geq & \Phi(u_n)\\
& \geq & \Phi(tu_n^+)\\
& = & \frac{1}{2}t^2\|u_n^+\|^2 + \frac{t^2}{2}\int_{\mathbb{R}^3}V(x)|u_n^+|^2dx - \int_{\mathbb{R}^3}K(x)F(t|u_n^+|)dx\\
& \geq & \frac{t^2\delta}{2} + o_n(1),
\end{eqnarray*}
which is a clear contradiction.

Then there exists $t_u > 0$ and $v_u \in E^-$ such that $t_u u^+ + v_u \in \mathcal{M}$. Note that by Proposition \ref{embedding} and Lebesgue Dominated Convergence Theorem, since $u_n \in \mathcal{M}$, it follows that
\begin{eqnarray*}
c & \leq & \Phi(t_u u^+ + v_u) \\
& = & \frac{t_u^2}{2}\|u^+\|^2 - \frac{1}{2}\|v_u\|^2 + \frac{1}{2}\int_{\mathbb{R}^3}V(x)|t_u u^+ + v_u|^2dx- \int_{\mathbb{R}^3}K(x)F(|t_u u^+ + v_u|)dx\\
& \leq & \liminf_{n \to \infty}\left(\frac{t_u^2}{2}\|u_n^+\|^2 - \frac{1}{2}\|v_u\|^2 + \frac{1}{2}\int_{\mathbb{R}^3}V(x)|t_u u_n^+ + v_u|^2dx\right.\\
& & \hspace{1.5cm}\left. - \int_{\mathbb{R}^3}K(x)F(|t_u u_n^+ + v_u|)dx\right)\\
& = & \liminf_{n \to \infty} \Phi(t_uu_n^+ + v_u)\\
& \leq & \Phi(u_n)\\
& = & c,
\end{eqnarray*}
which implies that the infimum of $\Phi$ on $\mathcal{M}$ is achieved in $t_u u^+ + v_u$ and finish the proof.
\vspace{1cm}

\noindent {\bf Acknowledgments:} This work has been finished while G. M. Figueiredo was as a Visiting Professor in FCT - Unesp. He would like to express his gratitude by the warm hospitality.


\begin{thebibliography}{99}

\bibitem{AlvesSouto} Alves, C.O., Souto, M.A., {\it Existence of solutions for a class of nonlinear
Schr\"{o}dinger equations with potential vanishing at infinity}, J.
Differential Equations, 254 (2) (2013), 1977-1991.

\bibitem{FigueiredoBarile} Barile, S., Figueiredo, G.M., {\it Existence of least energy positive, negative and nodal solutions for a class of $p\&q$-problems with potentials vanishing at infinity}, J. Mat. Anal. Appl., 427 (2) (2015), 1205 - 1233.

\bibitem{DingBartsch} Bartsch, T., Ding, Y., {\it Solutions of nonlinear Dirac equations}, J. Differential Equations, 226 (2006), 210 - 249.

\bibitem{DelPino} Del Pino, M., Felmer, P. {\it Local mountain pass for semilinear elliptic problems in unbounded domains}, Calc. Var. Partial
Differential Equations, 4 (1996), 121 -
137.

\bibitem {LivroDing} Ding, Y., {\it Variational methods for strongly indefinite problems}, Reviews in Mathematical Physics, 24 (10), (2012).

\bibitem {DingLiu} Ding, Y., Liu, X., {\it On semiclassical ground states of a nonlinear Dirac equation}, Interdisciplinary mathematical sciences, Vol - 7, World Scientific, (2007).

\bibitem {DingLiu2} Ding, Y., Liu, X., {\it Semi-classical limits of ground states of a nonlinear Dirac equation}, J. Differential Equations, 252, (2012), 4962 - 4987.

\bibitem {DingRuf} Ding, Y., Ruf, B., {\it Solutions of a nonlinear Dirac equation with external fields}, Arch. Rational Mech. Anal, 190 (2008), 1007 - 1032.

\bibitem {DingXu} Ding, Y., Xu, T., {\it Localized concentration of semi-classical states for nonlinear Dirac equations}, Arch. Rational Mech. Anal, 216 (2015), 415 - 447.

\bibitem{SereEsteban} Esteban, M., S\'er\'e, E., {\it Stationary states of the nonlinear Dirac equation: A variational approach}, Comm. Math. Phys. (1995), 323 - 350.

\bibitem{Merle} F. Merle, {\it Existence of stationary states for Dirac equations}, J. Differential Equations, 74, (1988), 50-68.

\bibitem{Nezza} Nezza, E., Palatucci, G., Valdinoci, E. {\it Hitchhiker's guide to the fractional Sobolev spaces}, Bull. Sci. Math., 136, (2012), 521 - 573.

\bibitem{Pankov} A. Pankov, {\it On decay of solutions to nonlinear Schr\"odinger equations}, Proc. Amer. Math.
Soc., 136, (2008), 2565?2570.

\bibitem{Rabinowitz} Rabinowitz, P. H., {\it On a class of nonlinear Schr\"{o}dinger equations}, ZAMP, 43, (1992), 270 - 291.

\bibitem{SzulkinWeth} Szulkin, A., Weth, T., {\it The method of Nehari manifold}, Handbook of nonconvex analysis and
applications, Int. Press., Somerville, (2010, 597 - 632.

\bibitem{Willem} Willem, M., {\it Minimax methods }, Handbook of nonconvex analysis and
applications, Int. Press., Somerville, (2010, 597 - 632.

\bibitem {Zhang1} Zhang, J., Tang, X., Zhang W., {\it Ground state solutions for nonperiodic Dirac equation with superquadratic nonlinearity},  J. Math. Phys., 54 (2013), 101502.

\bibitem {Zhang2} Zhang, J., Tang, X., Zhang W., {\it On ground state solutions for superlinear Dirac equation},  Acta Mathematica Scientia, 34B, (2014), 840 - 850.

\end{thebibliography}
\end{document}